\documentclass{tran-l}
\usepackage{amssymb,amsmath,latexsym}
\usepackage{cyrillic}

\newtheorem{thm}{Theorem}[section]
\newtheorem{cor}[thm]{Corollary}

\newtheorem{prop}[thm]{Proposition}
\theoremstyle{definition}

\theoremstyle{remark}
\newtheorem{rem}[thm]{Remark}
\numberwithin{equation}{section}

\newcommand{\cyrrm}{\fontencoding{OT2}\selectfont\textcyrup}

\begin{document}

\title[On the rank of a family of elliptic curves]
 {On the rank of certain parametrized elliptic curves }

\author{ A. Astaneh-Asl }

\address{Department of Mathematics, Arak-Branch, Islamic
Azad \\University, Arak, P.O. BOX 38135/567, Iran.}

\email{astaneh-asl@iau-arak.ac.ir; astanehasl@yahoo.com}

\thanks{Research supported by Arak-Branch Islamic Azad University.}

\subjclass{11G05.}

\keywords{ Elliptic Curve, Selmer Group.}


\dedicatory{}

\commby{Daniel J. Rudolph}


\begin{abstract}
In this paper the family of elliptic curves over $\mathbb{Q}$
given by the equation $ E_{p}:\ Y^2=(X-p)^3+X^3+(X+p)^3$ where $p$
is a prime number, is studied. It is shown that the maximal rank of the elliptic curves
is at most $3$ and some conditions under which we have  $\mbox{rank}(E_{p}(\mathbb{Q}))=0$
or  $\mbox{rank}(E_{p}(\mathbb{Q}))=1$ or  $\mbox{rank}(E_{p}(\mathbb{Q}))\geq2$ are given. Moreover
It is shown that in the family there are elliptic curves with rank 0,1,2 and 3.
\end{abstract}

\maketitle

\section*{Introduction}

Let $E$ be an elliptic curve over $\mathbb{Q}$ and $E(\mathbb{Q})$ be
the Mordell-Weil group of $E$ over $\mathbb{Q}$ which is finitely generated abelian group.
 The rank of $E(\mathbb{Q})$ as a $\mathbb{Z}$-module is called the rank of $E$ over
$\mathbb{Q}$. There is no algorithm which can compute the rank of any given elliptic curve so far.
 So it seems necessary to consider
certain families of elliptic curves and investigate their ranks
(see \cite{hsq,sc,st,wa}).
\par In this paper we consider the family of elliptic curves over $\mathbb{Q}$ given by the
equation
 $$E_{p}:\  Y^2=(X-p)^3+X^3+(X+p)^3,$$
where $p$ is a prime number, and prove the following theorem and propositions.

\begin{thm}\label{th1}
We have
\begin{enumerate}
\item If $p=2,3\ or\ p\equiv7\ (\mbox{mod}\ 24)$, then $\mbox{rank}(E_{p}(\mathbb{Q}))=0.$
\item If $p\equiv5,13,17\ (\mbox{mod}\ 24)$, then $\mbox{rank}(E_{p}(\mathbb{Q}))\leq1.$
\item If $p\equiv1\ (\mbox{mod}\ 24)\ \&\ (\frac{2}{p})_4=1$, then $\mbox{rank}(E_{p}(\mathbb{Q}))\leq3.$
\item In the other cases, $\mbox{rank}(E_{p}(\mathbb{Q}))\leq2.$
\end{enumerate}
\end{thm}

\begin{prop}\label{pr0}
Let $p\equiv17\ (\mbox{mod}\ 24)$ and $(\frac{2}{p})_4=1$. If there are integers $a, b$ such that
$3p=a^{4}+2b^{4}$, then $\mbox{rank}(E_{p}(\mathbb{Q}))=1.$
\end{prop}

\begin{prop}\label{pr00}
Let $p\equiv1\ (\mbox{mod}\ 24)$ and $(\frac{2}{p})_4=1$. If there are integers $a, b, c, \mbox{and}\ d$
such that $p=a^{4}+18b^{4}$ and
$3p=c^{4}+2d^{4}$, then $\mbox{rank}(E_{p}(\mathbb{Q}))\geq2.$
\end{prop}
Some primes which satisfy the conditions in Proposition \ref{pr0} and Proposition \ref{pr00}
will be given.


\section{The rank of $E_{p}$}

In this section Theorem \ref{th1} will be proved.
For proofing the theorem we have to deal with the Selmer groups
of $E_{p}$ corresponding to certain 2-isogenies. Let
$E/\mathbb{Q}$ be an elliptic curve over $\mathbb{Q}$ with a
torsion point of order 2. We use 2-descent via 2-isogeny method
 to compute the rank of $E$
over $\mathbb{Q}$ which is based on computing of the Selmer groups
corresponding to certain 2-isogeny of E (see \cite{cr,si,s-t}).

\par Let
$$E: y^{2}=x^{3}+ax^{2}+bx \ \ \  \ \ , a,b\in \mathbb{Z},$$
be an elliptic curve over $\mathbb{Q}$ and
$$\overline{E}: Y^{2}=X^{3}+\overline{a}X^{2}+\overline{b}X,$$
where $\overline{a}=-2a$ and $\overline{b}=a^{2}-4b$, be its
2-isogeny curve. Let

\begin{align}
 \Psi: E&\rightarrow\overline{E} \nonumber  \\
(x,y)&\longmapsto(\frac{y^{2}}{x^{2}}, \frac{y(b-x^{2})}{x^{2}}) \nonumber
\end{align}
be 2-isogeny of degree 2 and
\begin{align}
 \overline{\Psi}: \overline{E}&\rightarrow E \nonumber \\
(X,Y)&\longmapsto (\frac{Y^{2}}{4X^{2}}, \frac{Y(b-X^{2})}{8X^{2}}) \nonumber
\end{align}

be dual isogeny of $\Psi$. Let
$$C_{b_{1}}: b_{1}w^{2}=b_{1}^{2}+b_{1}\overline{a}z^{2}+bz^{4},$$
and
$$C_{\overline{b}_{1}}: \overline{b}_{1}w^{2}=
\overline{b}_{1}^{2}-2\overline{b}_{1}\overline{a}z^{2}+\overline{b}z^{4},$$
where $b_{1}|b$ and $\overline{b}_{1}|\overline{b}$, be the
homogeneous spaces for $E/\mathbb{Q}$ and
$\overline{E}/\mathbb{Q}$, respectively. The Selmer groups
corresponding to the 2-isogneies $\overline{\Psi}$ and $\Psi$ of
these curves are
$$S[\overline{\Psi}]=\{1.{{\mathbb{Q}}^{*}}^{2},
 b.{{\mathbb{Q}}^{*}}^{2}
\}\cup
 \{b_{1}.{{\mathbb{Q}}^{*}}^{2}:\ b_{1}|b\ and\
C_{b_{1}}(\mathbb{Q}_{p})\neq\phi\ \mbox{for all}\ p\in{S}\},$$

where $S:=\{\infty\}\cup\{p: \ p\ \mbox{is a prime and} \
p|2b\overline{b} \}$.\\ And
$$S[\Psi]=\{1.{{\mathbb{Q}}^{*}}^{2}, \overline{b}.{{\mathbb{Q}}^{*}}^{2} \}\cup
\{\overline{b}_{1}.{{\mathbb{Q}}^{*}}^{2}:\
\overline{b}_{1}|\overline{b}\ and\
C_{\overline{b}_{1}}(\mathbb{Q}_{p})\neq\phi\ \mbox{for all}\
p\in{S}\}.$$
\par
Now consider the elliptic curve $E_{p}$. With the change of variables $x=3X$ and $y=3Y$ the equation of
$E_{p}$ becomes
$$E_{p}: y^{2}=x^{3}+18p^{2}x,$$
and the following propositions
give us the structure of the Selmer groups.

\begin{prop}\label{pr1}
 Using the notations introduced above, we have
\begin{enumerate}
\item If $p\equiv11,19\ (\mbox{mod}\ 24)$ or $[p\equiv1\ (\mbox{mod}\ 24)\ \&\ (\frac{2}{p})_4=1]$,
 then $S_{p}[\overline{\Psi}]\cong(\frac{\mathbb{Z}}{2\mathbb{Z}})^{3}.$
\item If $p=2,3\ or\ p\equiv7\ (\mbox{mod}\ 24)$, then $S_{p}[\overline{\Psi}]\cong\frac{\mathbb{Z}}{2\mathbb{Z}}.$
\item In the other cases, $S_{p}[\overline{\Psi}]\cong(\frac{\mathbb{Z}}{2\mathbb{Z}})^{2}.$
\end{enumerate}
\end{prop}

\begin{prop}\label{pr2}
 We have
\begin{enumerate}
\item If $p\equiv23\ (\mbox{mod}\ 24)$ or $[p\equiv1\ (\mbox{mod}\ 24)\ \&\ (\frac{2}{p})_4=1]$,
 then $S_{p}[\Psi]\cong(\frac{\mathbb{Z}}{2\mathbb{Z}})^{2}.$
\item In the other cases, $S_{p}[\Psi]\cong\frac{\mathbb{Z}}{2\mathbb{Z}}.$
\end{enumerate}
\end{prop}

We proof Proposition \ref{pr1} and one can proof Proposition
\ref{pr2} by the same method.

By the definition it is clear that
$$\{1.{{\mathbb{Q}}^{*}}^{2},2.{{\mathbb{Q}}^{*}}^{2}\}\subseteq S_{p}[\overline{\Psi}].$$
So it is sufficient to check solvability of the equations
$$C_{b_{1}}: w^{2}=b_{1}+\frac{18p^{2}}{b_{1}}z^{4},$$
for $b_{1}=-1,\pm 3,\pm p,\pm 3p$ over $\mathbb{Q}_{l}$ where $l\in\{\infty,2,3,p\}$. If $b_{1}<0$ it is clear that
$C_{b_{1}}(\mathbb{Q}_{\infty})=\phi$ and then
 $b_{1}.{{\mathbb{Q}}^{*}}^{2}\notin S_{p}[\overline{\Psi}].$ For any $b_{1}>0$ we have
 $C_{b_{1}}(\mathbb{Q}_{\infty})\neq\phi.$ Let $p\neq 2,3$, we consider the equation

\begin{equation}\label{eq1}
 w^{2}=3+6p^{2}z^{4},
\end{equation}

corresponding to $b_{1}=3$. The solution $(z,w)=(1,1)$ for the
congruence $w^{2}\equiv 3+6p^{2}z^{4}$ (mod 8) can be lifted to a
solution for the equation (\ref{eq1}) in $\mathbb{Z}_{2}$ by using
Hensel's lemma, and then $C_{3}(\mathbb{Q}_{2})\neq\phi.$ By
considering the equation (\ref{eq1}) modulo 3 we have
$w^{2}\equiv 0$ (mod 3), say $w=3W$, so it can be written as
$3W^{2}=1+2p^{2}z^{4}$, and again the solution $(z,W)=(1,1)$ for
the congruence $3W^{2}\equiv 1+2p^{2}z^{4}$ (mod 3) lifts to a
solution for the equation $3W^{2}=1+2p^{2}z^{4}$ in $\mathbb{Z}_{3}$ which implies
$C_{3}({\mathbb{Q}_{3}})\neq\phi.$
 Now consider the equation (\ref{eq1}) over $\mathbb{Q}_{p}$.

\begin{enumerate}
 \item Let $(\frac{3}{p})=1$, then there is $w_{0}\in\mathbb{Z}$ such that
 $w_{0}^{2}\equiv3\ (\mbox{mod}\ p)$. The solution $(z,w)=(1,w_{0})$
 for the congruence $w^{2}\equiv 3+6p^{2}z^{4}$ (mod p) lifts to a solution for the
 equation (\ref{eq1}) in $\mathbb{Q}_{p}$.
 \item Let $(\frac{3}{p})=-1\ \&\ (\frac{2}{p})=-1$, then $(\frac{6}{p})=1$ and there is $w_{0}\in\mathbb{Z}$ such that
 $w_{0}^{2}\equiv6\ (\mbox{mod}\ p)$. Let $j$ be a positive integer number, then
 The solution $(z,w)=(1,w_{0})$ for the congruence $w^{2}\equiv 3p^{2+4j}+6z^{4}$ (mod p) lifts
 to a solution such as $(z,w)=(\alpha,\beta)$ for the equation $w^{2}=3p^{2+4j}+6z^{4}$ in $\mathbb{Q}_{p}$. So
 $(z,w)=(p^{-(1+j)}\alpha ,p^{-(1+2j)}\beta)$ is a solution for the equation (\ref{eq1}) in $\mathbb{Q}_{p}$.
 \item Let $(\frac{3}{p})=-1\ \&\ (\frac{2}{p})=1$, in this case one can show that there is no solution for the
 equation (\ref{eq1}) in $\mathbb{Q}_{p}$ since 3 and 6 are non-square mod p.
\end{enumerate}

Therefore
$$ p\equiv7,17\ (mod\ 24)\Leftrightarrow 3.{{\mathbb{Q}}^{*}}^{2}\notin
 S_{p}[\overline{\Psi}].$$

Now we deal with the case $b_{1}=p$. The corresponding equation is

\begin{equation}\label{eq2}
 w^{2}=p+18pz^{4}.
\end{equation}

Suppose that $C_{p}(\mathbb{Q}_{2})\neq\phi$. Since $v_{2}(w^{2})$ is even and $v_{2}(18pz^{4})$ is odd,
then necessarily $z,\ w\in\mathbb{Z}_{2}$, and therefore we deduce $p\equiv1,3\ (\mbox{mod}\ 8)$. Conversely
when $p\equiv1,3\ (\mbox{mod}\ 8)$ the solutions $(z,w)=(2,1)$ and $(z,w)=(1,1)$ for the congruence
$ w^{2}\equiv p+2pz^{4}$ (mod 8) lifts to  solutions for the equation (\ref{eq2}) in  $\mathbb{Q}_{2}$, respectively.

\par When $p\equiv1\ (\mbox{mod}\ 3)$ the solution $(z,w)=(1,1)$ for the equation (\ref{eq2}) mod 3 lifts to
a solution for it in $\mathbb{Q}_{3}$. Now let $p\equiv2\ (\mbox{mod}\ 3)$ and $j$ be a positive integer number,
the solution $(z,w)=(1,1)$ for the congruence $ w^{2}\equiv 3^{2+4j}p+2pz^{4}$ (mod 3) lifts to a solution
such as $(z,w)=(\alpha,\beta)$  for the equation  $ w^{2}=3^{2+4j}p+2pz^{4}$ in $\mathbb{Q}_{3}$.
So $(z,w)=(3^{-(1+j)}\alpha, 3^{-(1+2j)}\beta)$ is a solution for the equation (\ref{eq2}) in $\mathbb{Q}_{3}$, therefore
$C_{p}(\mathbb{Q}_{3})\neq\phi$.

\par Let $p\equiv3\ (\mbox{mod}\ 8)$, then there is $z_{0}\in\mathbb{Z}$ such that $1+18z_{0}^{4}\equiv0\ (\mbox{mod}\ p)$
since for any integer $x$, one of $x$ and $-x$ is a quadratic residue and the other one is a non-residue. So the solution
 $(z,W)=(z_{0},1)$ for the congruence $pW^{2}\equiv 1+18z^{4}\ (\mbox{mod}\ p)$ lifts to a solution for the equation
$pW^{2}=1+18z^{4}$ in $\mathbb{Q}_{p}$ and then $C_{p}(\mathbb{Q}_{p})\neq\phi$. Now let $p\equiv1\ (\mbox{mod}\ 8)$
and $C_{p}(\mathbb{Q}_{p})\neq\phi$. Suppose that $(z,w)$ is a solution for the equation (\ref{eq2}) in $\mathbb{Q}_{p}$.
Let $v_{p}(w)=k$ and  $v_{p}(z)=j$, it is clear that $j$ must be zero and $k>0$. So considering equation (\ref{eq2}) mod p
implies that $(\frac{-18}{p})_{4}=1$ where $(\frac{}{p})_{4}$ is the rational quartic residue symbol mod p. Conversely, if
 $(\frac{-18}{p})_{4}=1$ then $C_{p}(\mathbb{Q}_{p})\neq\phi$. Therefore $p.{\mathbb{Q}^{*}}^{2}\in S_{p}[\overline{\Psi}]$
 if and only if
 $$[p\equiv11,19\ (\mbox{mod}\ 24)]\ or\ [p\equiv1\ (\mbox{mod}\ 24)\ \&\ (\frac{2}{p})_{4}=1 ]\ or\
  [p\equiv17\ (\mbox{mod}\ 24)\ \&\ (\frac{2}{p})_{4}=-1].$$
In the case $b_{1}=3p$, the corresponding equation is

\begin{equation}\label{eq3}
 w^{2}=3p+6pz^{4}.
\end{equation}
By the same methods as in the case  $b_{1}=p$ one can show that
$$[p\equiv11,19\ (\mbox{mod}\ 24)]\ or\ [p\equiv1,17\ (\mbox{mod}\ 24)\ \&\ (\frac{2}{p})_{4}=1 ]\\
\Leftrightarrow 3p.{\mathbb{Q}^{*}}^{2}\in S_{p}[\overline{\Psi}].$$
Finally for $p=2,3$ we have
$$S_{2}[\overline{\Psi}]=S_{3}[\overline{\Psi}]=\{1.{{\mathbb{Q}}^{*}}^{2},2.{{\mathbb{Q}}^{*}}^{2} \},$$
which completes the proof of Proposition \ref{pr1}. $\ \ \Box$

\begin{cor}\label{col}
We have
\begin{enumerate}
\item If $p\equiv11,19\ (\mbox{mod}\ 24)$ or $[p\equiv1\ (\mbox{mod}\ 24)\ \&\ (\frac{2}{p})_4=1]$,
 then $S_{p}[\overline{\Psi}]=\{ b_{1}.{{\mathbb{Q}}^{*}}^{2}:\ b_{1}=1, 2, 3, 6, p, 2p, 3p, 6p\}.$
\item If $[p\equiv5,13,23\ (\mbox{mod}\ 24)]\ or\ [p\equiv1\ (\mbox{mod}\ 24)\ \&\ (\frac{2}{p})_{4}=-1 ]$, then $S_{p}[\overline{\Psi}]=\{ b_{1}.{{\mathbb{Q}}^{*}}^{2}:\ b_{1}=1, 2, 3, 6\}.$
\item If $p\equiv17\ (\mbox{mod}\ 24)\ \&\ (\frac{2}{p})_{4}=1 $, then $S_{p}[\overline{\Psi}]= \{ b_{1}.{{\mathbb{Q}}^{*}}^{2}:\ b_{1}=1, 2, 3p, 6p\}.$
 \item If $p\equiv17\ (\mbox{mod}\ 24)\ \&\ (\frac{2}{p})_{4}=-1 $, then $S_{p}[\overline{\Psi}]=\{ b_{1}.{{\mathbb{Q}}^{*}}^{2}:\ b_{1}=1, 2, p, 2p\}.$
\item In the other cases, $S_{p}[\overline{\Psi}]=\{1.{{\mathbb{Q}}^{*}}^{2},2.{{\mathbb{Q}}^{*}}^{2} \}.$
\end{enumerate}
And\\

$
 S_{p}[\Psi]=\left\{
  \begin{array}{ll}
    \{ b_{1}.{{\mathbb{Q}}^{*}}^{2}:\ b_{1}=1, -2, p, -2p\}\ \ , & \hbox{$ p\equiv1\ (\mbox{mod}\ 24)\ \&\ (\frac{2}{p})_{4}=1 $;}\\
    \{ b_{1}.{{\mathbb{Q}}^{*}}^{2}:\ b_{1}=1, -2, -p, 2p\}\ \ , & \hbox{$p\equiv23\ (\mbox{mod}\ 24)$;} \\
    \{1.{{\mathbb{Q}}^{*}}^{2},-2.{{\mathbb{Q}}^{*}}^{2} \}, & \hbox{$\mbox{otherwise}$.}
  \end{array}
\right.
$

\end{cor}

Note that the first result in Corollary \ref{col} is clear because of the proof of Proposition \ref{pr1}, and one can obtain
the second part with the same method.

Now we have the structures of the Selmer groups and we can proof
 Theorem \ref{th1}. Consider the following map

\begin{align}
\alpha_{p}:\ E_{p}(\mathbb{Q})&\rightarrow S_{p}[\overline{\Psi}] \nonumber\\
 \mathcal{O} &\longmapsto 1.{{\mathbb{Q}}^{*}}^{2} \nonumber\\
(0,0) &\mapsto 2.{{\mathbb{Q}}^{*}}^{2} \nonumber\\
(x,y) &\mapsto x.{{\mathbb{Q}}^{*}}^{2}\ \ \ for\ x\neq0.\nonumber
\end{align}

The following sequence is exact
$$0\rightarrow E_{p}(\mathbb{Q})/\overline{\Psi}(\overline{E}_{p}(\mathbb{Q}))\rightarrow
 S_{p}[\overline{\Psi}]\rightarrow {\cyrrm{SH}}_{p}[\overline{\Psi}]\rightarrow0 $$
 where ${\cyrrm{SH}}_{p}[\overline{\Psi}]$ is the cokernel of the left hand
 side injection which is called Tate-Shafarevich group of $E_{p}$.
 For the rank of $E_{p}$ and $\overline{E}_{p}$ one obtains the following
 formula
 $$\mbox{rank}(E_{p}(\mathbb{Q}))=dim_{\mathbb{F}_{2}}(S_{p}[\overline{\Psi}])+
 dim_{\mathbb{F}_{2}}(S_{p}[\Psi])-dim_{\mathbb{F}_{2}}({\cyrrm{SH}}_{p}[\overline{\Psi}])
 -dim_{\mathbb{F}_{2}}({\cyrrm{SH}}_{p}[\Psi])-2.$$
Now one can easily complete the proof of Theorem \ref{th1}.  $\ \ \ \Box$
\par
For proofing Proposition \ref{pr0} and Proposition \ref{pr00} we use the same method as in \cite{k-m}, for more details see
\cite{ch} or \cite{s-t}.\\

\textbf{Proof of Proposition \ref{pr0}}.
Since $p\equiv17\ (\mbox{mod}\ 24)\ \mbox{and}\ (\frac{2}{p})_{4}=1$ by Corollary \ref{col} we have
$S_{p}[\overline{\Psi}]=\{ 1.{{\mathbb{Q}}^{*}}^{2},2.{{\mathbb{Q}}^{*}}^{2},3p.{{\mathbb{Q}}^{*}}^{2},6p.{{\mathbb{Q}}^{*}}^{2}\}$
and $S_{p}[\Psi]=\{ 1.{{\mathbb{Q}}^{*}}^{2},-2.{{\mathbb{Q}}^{*}}^{2}\}.$
Consider the equation $3pS^{4}+6pT^{4}=U^{2}$, $(S,T,U)=(a,b,3p)$ is a solution for it such that $a,\ b\geq1$ and $gcd(a,6p)=1$.
Because $a$ is odd and if $gcd(a,6p)=d$, then $d$ will be odd. Now note that $6p=2a^{4}+4b^{4}$  so $d|4b^{4}$ and then $d|b^{4}$
therefore we deduce $d|p$. Suppose $d=p$, this implies that $d|a$ and $d|b$, i.e., there are integers $a_{1}$ and $b_{1}$ such that
$a=pa_{1}$ and $b=pb_{1}$, so we have $3p=p^{4}(a_{1}^{4}+2b_{1}^{4})$ which is a contradiction. Thus $d=1$ and then
$3p.{\mathbb{Q}^{*}}^{2}\in Im\alpha_{p}$.
Therefore $Im\alpha_{p}=\{1.{\mathbb{Q}^{*}}^{2},2.{\mathbb{Q}^{*}}^{2},3p.{\mathbb{Q}^{*}}^{2},6p.{\mathbb{Q}^{*}}^{2}\}$ and
$Im\overline{\alpha}_{p}=\{1.{\mathbb{Q}^{*}}^{2},-2.{\mathbb{Q}^{*}}^{2}\}$, and then $rank(E_{p}(\mathbb{Q}))=1$.$\ \ \Box$
\vspace{4mm}
\par Note that $p=1217,\ 1601,\ 5297,\ 9521$ are some primes which satisfy the conditions in Proposition \ref{pr0}.
\vspace{4mm}

\textbf{Proof of Proposition \ref{pr00}}.
By  Corollary \ref{col} we have
$$S_{p}[\overline{\Psi}]=\{ 1.{{\mathbb{Q}}^{*}}^{2},2.{{\mathbb{Q}}^{*}}^{2}, 3.{{\mathbb{Q}}^{*}}^{2}, 6.{{\mathbb{Q}}^{*}}^{2}, p.{{\mathbb{Q}}^{*}}^{2},2p.{{\mathbb{Q}}^{*}}^{2},3p.{{\mathbb{Q}}^{*}}^{2},6p.{{\mathbb{Q}}^{*}}^{2}\}$$
and $S_{p}[\Psi]=\{ 1.{{\mathbb{Q}}^{*}}^{2},-2.{{\mathbb{Q}}^{*}}^{2}, p.{{\mathbb{Q}}^{*}}^{2}, -2p.{{\mathbb{Q}}^{*}}^{2}\}$, and
Proposition \ref{pr0} implies that
$\{1.{\mathbb{Q}^{*}}^{2},2.{\mathbb{Q}^{*}}^{2},3p.{\mathbb{Q}^{*}}^{2}, 6p.{\mathbb{Q}^{*}}^{2}\}\subseteq Im{\alpha}_{p}$.
On the other hand $(S,T,U)=(a,b,p)$ is a solution for the equation $pS^{4}+18pT^{4}=U^{2}$ which implies
that $p.{\mathbb{Q}^{*}}^{2}\in Im\alpha_{p}$, so $\#Im\alpha_{p}=\#S_{p}[\overline{\Psi}]=8$ and then $rank(E_{p}(\mathbb{Q}))\geq2.$
Note that $gcd(a,18p)=1$, since if $gcd(a,18p)=d_{1}$, then $d_{1}|p$ because $a$ is odd and $3\nmid a$. Suppose $d_{1}=p,$
this concludes that $d_{1}|a$ and $d_{1}|b$ which give us a contradiction. $\ \ \Box$

\vspace{4mm}
 $p=19249$ is a prime which satisfies the conditions in Proposition \ref{pr00}.
\begin{rem}
By using Sage one can see that $rank(E_{7}(\mathbb{Q}))=0$, $rank(E_{5}(\mathbb{Q}))=1$,
 $rank(E_{11}(\mathbb{Q}))=2$ and $rank(E_{19249}(\mathbb{Q}))=3$.
\end{rem}


\end{document}